\documentclass[11pt]{article}
\usepackage{amsfonts, amsmath,graphicx, amsthm, amssymb, cite, stmaryrd,mathtools}
\numberwithin{equation}{section}

\numberwithin{proposition}{section}
\newtheorem{thm}{Theorem}
\numberwithin{thm}{section}

\numberwithin{definition}{section}
\newtheorem{rem}{Remark}
\numberwithin{rem}{section}
\newtheorem{lem}{Lemma}
\numberwithin{lem}{section}

\numberwithin{cor}{section}
\newcommand{\RR}{\mathrm{I\!R\!}}
\newcommand{\TT}{\mathrm{I\!T\!}}
\usepackage[margin=1in]{geometry}
\usepackage{color}

\usepackage{pgfplots}

\usetikzlibrary{arrows.meta}
\usetikzlibrary{shapes,backgrounds}
\usetikzlibrary{patterns}

\pgfplotsset{compat = newest}

\pgfplotsset{
  every axis/.append style={
    axis x line=middle,    
    axis y line=middle,    
    axis line style={<->,color=blue}, 
    xlabel={$x$},          
    ylabel={$y$},          
  }
}

\begin{document}

\title{Non-Uniqueness in Plane Fluid Flows}
\author{Heiko Gimperlein\thanks{Engineering Mathematics, University of Innsbruck, Innsbruck, Austria} \thanks{Department of Mathematical, Physical and Computer Sciences, University of Parma,
43124 Parma, Italy} \and Michael Grinfeld\thanks{Department of Mathematics and Statistics, University of Strathclyde, Glasgow, G1 1XH, UK} \and Robin J.~Knops\thanks{The Maxwell Institute of Mathematical Sciences and School of Mathematical and Computing Sciences, Heriot-Watt University, Edinburgh, EH14 4AS, Scotland, UK} \and Marshall Slemrod\thanks{Department of Mathematics, University of Wisconsin, Madison, WI 53706, USA}}
\date{}
\maketitle

\providecommand{\keywords}[1]{{\textit{Key words:}} #1}
\providecommand{\msc}[1]{{\textit{MSC classes:}} #1}

\begin{abstract}

\noindent Examples of dynamical systems proposed by Z.~Artstein 
and C.~M.~Dafermos 
admit non-unique solutions that track a one parameter family of closed circular orbits contiguous at a single point. Switching between orbits at this single point produces an infinite number of solutions with the same initial data. Dafermos appeals to a maximal entropy rate criterion to recover uniqueness. 

\noindent These results are here interpreted as non-unique Lagrange trajectories on a particular spatial region. The corresponding velocity is  proved consistent with plane steady compressible fluid flows that  for specified pressure and mass density  satisfy not only the Euler equations but also the  Navier-Stokes equations for specially chosen volume and   (positive) shear  viscosities. The maximal entropy rate criterion recovers uniqueness. 
\end{abstract}
\\
\msc{76N10 (primary), 34A12, 35Q35 (secondary).}

\section{Introduction}\label{intro}
This paper  derives 
  explicit non-unique continuous  Lagrange trajectories 
 related to  certain  Euler and Navier-Stokes steady compressible plane  fluid flow. 
 The orbits tracked by the trajectories belong to a one parameter family of closed contiguous circular paths. Non-uniqueness occurs when trajectories switch orbits at the common point of intersection.
Even so,  uniqueness follows for the problems considered here by appeal to an entropy rate criterion that isolates a single preferred trajectory.

 The non-unique behaviour is not dissimilar to that occurring in the convex integration investigations    by De Lellis and Sz\`{e}kelyhidi \cite{ds09, ds10} and Buckmaster and Vicol \cite{bv19} for the incompressible Euler and Navier-Stokes systems. 
 It is therefore reasonable to speculate  whether a unique solution also might  be selected  by application of an entropy rate criterion. 
 Brief comment on this aspect is given  in Section~\ref{ci}. 

The  present analysis is partly inspired  by  an example in control theory  proposed  by Artstein \cite[see eqn.~(6.1)]{a83} who, rather than non-uniqueness, emphasises  non-existence of suitable controls that  otherwise stabilise the system and ensure solutions asymptotically approach  zero. The primary  motivation, however, comes from  what Dafermos \cite{d12} refers to as an artificial example of a nonlinear oscillator governed by an ordinary differential equation.  The  associated  trajectories  considered by both Artstein and Dafermos are circular, and  touch  at their  common point of intersection taken to be the origin.
Consequently, trajectories  can be switched at the origin by instantaneous variation of an  \emph{entropy}   function that is constant along individual trajectories. In this manner, irrespective of initial conditions, an uncountable set of different solutions can be chosen to establish non-uniqueness. 


Dafermos \cite{d12} recovers  uniqueness   by selecting that particular solution which   dissipates  entropy   at  maximum rate. Such a  solution  is not only  unique but  physically relevant              in the sense 
of Hadamard's definition of well-posedness.   Similar entropy rate criteria, originally formulated to treat shock wave propagation, include 
the maximal entropy rate  proposed by  Ziegler \cite{z83} and the minimum entropy rate  due to Prigogine \cite[Chapt.~V]{p47}.             Both principles are  reviewed in the monograph by Moroz \cite[Chapt.~1]{m11}.  Glimm, Lazarev and Chen \cite{glc20} remark that the  subtle difference between these two principles   may cause confusion. Ziegler's principle is relevant to closed thermodynamic systems while Prigogine's principle applies to open thermodynamic systems. Evidently, Dafermos's  investigation relates to  a closed system. 

 Dafermos's example, however,  is far from trivial.   It is shown in this paper that  solutions,  interpreted as steady Lagrange  trajectories (particle paths),  are related to  steady  plane flows for compressible  Euler and Navier-Stokes equations with variable mass density and positive shear viscosity. The extension   to fluid problems  
demonstrates that  tracking arbitrarily prescribed discrete entropy 
profiles produces  an infinite number of continuous Lagrange trajectories  with the same initial conditions. Consequently, trajectories are not unique. 
Nevertheless, a single unique trajectory  is identified by the maximal entropy rate admissibility criterion proposed by Dafermos \cite{d12}. 
 An interpretation of this result is  that ``wild'' initial data can be  produced for steady solutions to an initial boundary problem for both  the Euler and Navier-Stokes equations. The term ``wild'' is  understood in the sense that there are  an infinite number of Lagrange trajectories for the same initial conditions. The entropy rate criterion identifies a unique trajectory and renders  the notion of ``wildness''  redundant. 
 
The conclusion contrasts with the uniqueness results for weak solutions to  the incompressible Navier-Stokes equations  obtained by 
Robinson and Sadowski \cite{rs09a, rs09b}. These authors show for sufficiently regular velocity that the Lagrange trajectories are unique and smooth in time for almost every initial data. On the other hand, the velocity for the  Lagrange trajectories considered here
lacks Lipschitz continuity at one point and therefore  does not   satisfy the continuity conditions imposed  by Robinson and Sadowski.  

 An alternative  admissibility  procedure, introduced by Brenier \cite{bb18},  relies upon the concave maximization of a certain functional to derive smooth solutions to  Euler's equation.  Further  procedures  include that by Lasarzik \cite{l22} who adopts  a maximal entropy rate to conclude that ``dissipative solutions'' as defined by P.-L. Lions \cite{pll96}  exist, are unique,  and depend continuously  upon initial data for both the Euler and Navier-Stokes equations of incompressible fluid dynamics. It must be noted, however, that such dissipative solutions are not in general weak solutions.  A different notion of dissipative solution is treated by Breit, Feireisl and Hofmanova \cite{bfh20} who replace the concept of a weak solution by a semi-flow satisfying a maximal entropy rate production criterion. By construction, the semi-flow implies modified well-posedness of the Euler system. 

Also of importance is whether the maximal entropy rate  criterion can be accurately simulated by numerical analysis.
Glimm, Lazarev and Chen \cite{glc20}  note that this  comparatively delicate issue is   of substantial physical importance. A  numerical method that preserves qualitative behaviour is required and, in fact, it is shown in Section~\ref{numapp}  that the symplectic Euler method successfully ensures convergence to the unique solution 
selected by the criterion.  The explicit Euler method fails in this respect.


 As already remarked,  based upon  the Artstein-Dafermos examples, the non-uniqueness   described in this paper  for the Euler and Navier-Stokes equations is due to the corresponding Lagrange trajectories switching at the origin from one circular path to another along each of which   there is  different constant \emph{entropy}.  Arbitrary prescription of these entropies  and the order in which the trajectories are described results in an uncountable number of distinct solutions satisfying the same initial conditions. The behaviour has   features in common with that respectively established by  De Lellis and Sz\`{e}kelyhidi \cite{ds09, ds10} and by  Buckmaster and Vicol \cite{bv19} for weak solutions on a torus to the incompressible Euler and Navier-Stokes equations. These authors use convex integration techniques to prove that  an arbitrary number of such weak solutions exist possessing both the same assigned smooth global energy   and the same initial conditions. Consequently, non-uniqueness is established. (Extension to the compressible Euler equations is due to Chiodaroli and Kreml \cite {ck14}, and to Feireisl \cite{ef14}).



Arbitrary  specification of an appropriate   profile is thus common to both developments and  invites exploration of the possible application of  a maximal entropy rate criterion  to identify a unique solution  from among the many  weak solutions determined by convex integration methods.  

 Of separate interest is the physical relevance of  weak solutions obtained by convex integration  to the Euler and Navier-Stokes equations. Such solutions are the limit of increasingly high frequency oscillatory perturbations that  eventually may contravene  the basic continuum hypotheses assumed in the  derivation of  the equations.  In this respect,  the Lagrange trajectories considered here  are explicit and indeed are continuously differentiable except at the origin.  Clearly, apart from this singularity, they accord with  continuum hypotheses.  



Section~\ref{mot} commences with   Dafermos's construction of non-unique orbits for a nonlinear oscillator that involves the arbitrary prescription of an \emph{entropy} function constant on each path. A theorem conveniently summarises relevant conclusions.  Two admissibility criteria proposed by Dafermos are then described that  recover a physically meaningful unique trajectory. The first criterion augments the system of equations  by a friction term. Solutions to the penalised system tend to a unique solution in the limit as friction tends to zero.   The second, or maximal entropy rate,  criterion requires maximum dissipation of the individual entropies  as  trajectories are successively switched and leads to the same path obtained by the limiting friction argument.  Section~\ref{ens} establishes that  the velocity occurring in  the nonlinear oscillator satisfies  the plane compressible Euler equations with  specified pressure and  particular non-uniform mass density. Moreover, the corresponding Lagrange trajectories are geometrically similar to those obtained from the nonlinear oscillator and therefore possess  similar non-uniqueness features.  Isolation of a unique trajectory follows from the maximal entropy rate criterion. The  analysis is extended to the corresponding plane compressible Navier-Stokes equations with specially chosen (positive) shear  and volume viscosities.  These  viscosity coefficients are computed in the Appendices. Section~\ref{numapp} discusses the numerical approximations of solutions  to the Dafermos system considered in Section~\ref{mot}.  The objective is to employ  a numerical method that not only preserves the qualitative behaviour of  constant entropy  along  individual trajectories, but   is also capable of simulating   trajectories  selected by  the maximal entropy rate criterion.   Convergence to a unique solution would then  be implied. The symplectic Euler method satisfies these requirements since it determines the approximation to the entropy function  to within an error of 
the order of the step size $h>0$. The approach is diagramatically illustrated and arguments are given that the discrete flow selects the solution given by the entropy rate admissibility criterion. Section~\ref{ci} comments on  the resemblance between the results of Section~\ref{ens} and those in the convex integration literature that  for a given smooth global energy construct an infinite number of solutions to the Euler equations on a torus. 
A final brief remark concerns the fundamental issue of whether the eventual irregularity of such solutions 
contravenes basic continuum hypotheses.


Standard notation is employed throughout the main text apart from Remark \ref{remark31} where bold type indicates vectors.  The Appendices introduce an indicial notation accompanied by the summation and comma conventions.


\section{Motivation and admissibility criteria}\label{mot}

This section reviews   previous contributions that inspired the present generalisation  to the Euler and Navier-Stokes equations. The discussion also explains how an  entropy rate criterion  recovers  a unique solution.

Let $(x(t),y(t))\in \RR^{2},\, t\in \RR,$ be state variables corresponding to the Cartesian coordinates of a point moving in the plane  as time varies.  Dafermos \cite{d12}   studies  motions whose state  variables satisfy the system of ordinary differential equations for a nonlinear oscillator  given by
\begin{eqnarray}
\label{orb1}
\dot{x}&=& y,\\
\label{orb3}
\dot{y}&=&\frac{(y^{2}-x^{2})}{2x},
\end{eqnarray} 
where a superposed dot indicates differentiation with respect to  time $t$.

Initial data are specified by
\begin{equation}
\label{iconds}
(x(0),\,y(0))=(x_{0},\,y_{0}),
\end{equation}
where $(x_{0},\,y_{0})$ are prescribed.

The velocity component \eqref{orb1} is continuously differentiable  everywhere  but the component \eqref{orb3} is not continuous at the origin. Dynamics are best understood from an integrated form of the equations given below.  

It is immediate  from \eqref{orb1} and \eqref{orb3}  that the entropy function $H(x,y)$, defined by
\begin{equation}
\label{hpprel}
H(x,y):=\frac{(x^{2}+y^{2})}{2x}, 
\end{equation}
 is invariant with respect to time $t$; that is 
\begin{equation}
\label{hder}
\frac{dH(x(t),y(t))}{dt}=0,
\end{equation}
or 
\begin{equation}
\label{hprel}
H(x,\,y)=c
\end{equation}
for  positive constant $c$.

Note that $H(x,y)$ is not to be confused with the (kinetic) energy. The relationship between these quantities is presented in Section~\ref{ci}.

The constant entropy function \eqref{hpprel} implies that 
the state variables traverse the one parameter family of closed circular orbits
\begin{equation}
\label{altorb}
x^{2}+y^{2}=2xc.
\end{equation}
Consequently, 
\begin{equation}
\label{orb}
(x-c)^{2}+y^{2}=c^{2},
\end{equation}
and consequently $0\le x(t)\le 2c,\, -c\le y(t)\le c.$ Moreover, \eqref{orb1} implies that $x(t)$ increases or decreases as $y(t)$ is positive or negative so that the closed circular orbits are traversed clockwise.  For varying  $c\ge 1$ they touch each other at their common point of intersection located at the origin.  Since solutions trace the family of orbits \eqref{orb},  equation \eqref{orb3} may be alternatively  written as
\begin{equation}
\label{orb2}
\dot{y}=(c-x).
\end{equation}

%


  Artstein's study  \cite[eqn.~6.1]{a83} is within the broader context of a control problem in  which  \eqref{orb1} and \eqref{orb2} are generalised to  the family 
\begin{eqnarray}
\label{a1}
\dot{x}&=&(2xy)w,\\
\label{a2}
\dot{y}&=&(y^{2}-x^{2})w,
\end{eqnarray}   
where  $w(x,y)$ is a scalar control  chosen to ensure stabilisability of the equilibrium state $x=y=0.$  It follows immediately that since
\begin{equation}
\label{aorb}
\frac{dx}{2xy}=\frac{dy}{(y^{2}-x^{2})},
\end{equation}
the paths followed by solutions to \eqref{a1} and \eqref{a2} are given by the family of circular  closed  orbits \eqref{orb} or alternatively \eqref{altorb}, and consequently the function $H(x,t)$ is constant along each path. Dependent upon the choice of $w$, the Lipschitz condition may fail 
which would imply  lack of uniqueness for   paths with given initial data and which we now examine.  




Particular choices of the control function $w(x,y)$ deserve attention.

First, set
\begin{equation}
\label{wgen}
w= x^{-\alpha},\qquad \alpha>0.
\end{equation}
Substitution in \eqref{a1} and \eqref{a2}  after appeal to \eqref{altorb} 
yields
\begin{eqnarray}
\label{genw1}
\dot{x}&=& 2x^{(\frac{3}{2}-\alpha)}(2c-x)^{1/2},\\
\label{genw2}
\dot{y}&=& 2x^{(1-\alpha)}(c-x).
\end{eqnarray}


Hence \eqref{altorb}, \eqref{genw1} represent an integrated form of \eqref{a1}, \eqref{genw1} where conservation of $H$ given \eqref{hprel}  is extended by continuity to the origin. It is here that lack of Lipschitz continuity of the right hand side of \eqref{genw1} makes itself apparent and the usual loss of uniqueness occurs in $x$ for $1\le \alpha < 3/2$ (see, for example, Burkill \cite{b75}), while $y$ is well defined from the choice of $x$. As stated below, $\alpha=1$ is the value adopted  by Dafermos \cite{d12}.


As second choice, put
\begin{equation}
\label{ah}
w=\frac{1}{2x^{2}},
\end{equation}
to obtain from  \eqref{a1} and \eqref{a2}  the Hamiltonian system
\begin{eqnarray}
\label{a3}
\dot{x}&=& \frac{y}{x},\\
\label{a4}
\dot{y}&=& \frac{(y^{2}-x^{2})}{2x^{2}},
\end{eqnarray}
whose Hamiltonian is the function  $H(x,y)$.

Finally,  consider 
\begin{equation}
\label{daf}
w=\frac{1}{2x}
\end{equation}
which reduces \eqref{a1} and \eqref{a2} to the system \eqref{orb1} and \eqref{orb3} considered by Dafermos \cite{d12}. Explicit solutions, of interest when interpreted in Section~\ref{ens} as Lagrange trajectories,  are easily derived. Details are presented in Appendix~\ref{orbcal}.  Express \eqref{orb} as
\begin{equation}
\nonumber
\dot{x}(t)=y(t)=\pm\sqrt{x(2c_{0}-x)},
\end{equation}
where $c$ has assumed the particular value $c_{0}$ specified from initial data according to 
\begin{equation}
\label{codef}
c_{0}=\frac{(x_{0}^{2}+y_{0}^{2})}{2x_{0}}.
\end{equation}
It is supposed that initial data are such that $c_{0}\ge 1.$  Choose the positive square root to obtain by integration 
\begin{equation}
\label{orbdafx}
x(t)-c_{0}=  c_{0}\cos{(-t+\theta_{0})}.
\end{equation}

The  corresponding expression for $y(t)$, obtained by differentiation of the last relation, becomes
\begin{equation}
\label{orbdafy}
y(t)=  c_{0}\sin{(-t+\theta_{0})},
\end{equation}
and shows that the constant $\theta_{0}$ is determined from initial data to be
\begin{equation}
\label{thdef}
\tan{\theta_{0}}=\frac{y_{0}}{x_{0}-c_{0}}.
\end{equation}

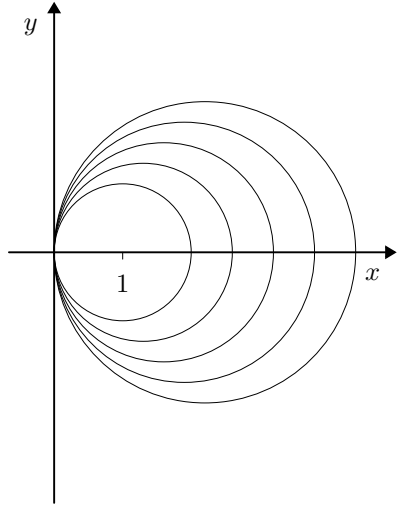
\begin{figure}[htb]
\centerline{
\resizebox{6cm}{!}{
\def\firstcircle{(1,0) circle (1.0cm)}
\def\secondcircle{(1.3,0) circle (1.3cm)}
\def\thirdcircle{(1.6,0) circle (1.6cm)}
\def\fourthcircle{(1.9,0) circle (1.9cm)}
\def\fifthcircle{(2.2,0) circle (2.2cm)}
\begin{tikzpicture}[
  line cap=round,
  line join=round,
  >=Triangle,
  myaxis/.style={->,thick}]

\begin{scope}
        \draw \firstcircle;
        \draw \secondcircle;
        \draw \thirdcircle;
        \draw \fourthcircle;
        \draw \fifthcircle;
    \end{scope}

\draw[myaxis] (-0.66,0) -- (5.00,0) node[below left = 1mm] {$x$};      
\draw[myaxis] (0,-3.66) -- (0,3.66) node[below left = 1mm] {$y$};

\draw (1,0)--(1,-.1) node[below = 1mm] {$1$};
\end{tikzpicture}
}
}
\caption{\label{phaseportrait}Phase portrait for \eqref{orbdafx} and \eqref{orbdafy}.}
\end{figure}

The circular orbits  \eqref{orbdafx} and \eqref{orbdafy} may otherwise be  derived by integration of \eqref{orb1} and \eqref{orb2}. They are centred at $(c_{0},0)$, where $c_{0}$ is determined from initial conditions \eqref{codef}, and are restricted to the set parametrized by $1\le c_{0}\le R$ sketched in Figure \ref{phaseportrait}.

Orbits tracked by the state variables,  as predicted by \eqref{orb}, are centred at $(c_{0},\,0)$ and pass through the origin  at times $(-t+\theta_{0}) =-(2n+1)\pi,\, n=0,1,2,\ldots$ . Each circular orbit is completed in the same time $2\pi$ but at speeds dependent upon  the radius:
\begin{equation}
\nonumber
\left(\dot{x}^{2}+\dot{y}^{2}\right)^{1/2}=c_{0},
\end{equation}
corresponding to rigid body rotation.

 Further properties are derived  upon conversion to polar coordinates $(r,\theta)$.   Set $x=r\cos{\theta},\,y=r\sin{\theta}$  so that  \eqref{orb1} and \eqref{orb3} become
\begin{eqnarray}
\label{dangde}
\dot{\theta}&=& -\frac{1}{2},\\
\label{dradde}
\dot{r}&=& \frac{1}{2} r\tan{\theta},
\end{eqnarray}
which by integration  lead to
\begin{eqnarray}
\label{dang}
\theta(t)&=& -\frac{t}{2}+\theta_{0}/2,\\
\label{drad}
r(t)\cos{(\theta_{0}/2)}&=& r_{0}\cos{\theta(t)},
\end{eqnarray}
where $\theta_{0}$ satisfies  \eqref{thdef} and $(r_{0},\,\theta_{0}/2)$, the initial values of $(r(t),\,\theta(t))$, accordingly are given by 
\begin{equation}
\nonumber
r_{0}=\left(x^{2}_{0}+y^{2}_{0}\right)^{1/2},\qquad \tan{(\theta_{0}/2)}=\frac{y_{0}}{x_{0}},
\end{equation}
and 
\begin{equation}
\nonumber
r_{0}=2c_{0}\cos{(\theta_{0}/2)}.
\end{equation}
The last expression inserted into \eqref{drad} gives
\begin{equation}
\label{rexp}
r(t)=2c_{0}\cos{\theta(t)},
\end{equation}
which is the circular orbit previously obtained.

The state variables  expressed by either \eqref{orbdafx} and \eqref{orbdafy}, or by \eqref{dang} and \eqref{rexp}, indicate the explicit  manner in which orbits pass through the origin $(0,0)$.

The calculations so far suppose that $c\geq 1$.   To include the region  $c< 1$,   Dafermos \cite{d12} considers the modified system
\begin{eqnarray}
\label{dorb1}
\dot{x}&=& y,\qquad \dot{y}=\frac{(y^{2}-x^{2})}{2x},\quad \mbox{for}\qquad (x-1)^{2}+y^{2}\geq 1,\\
\label{dorb2}
\dot{x}&=& y,\qquad \dot{y}= (1-x), \qquad (x-1)^{2}+y^{2}<1.
\end{eqnarray} 
Solutions to \eqref{dorb2}  are of the same form as \eqref{orbdafx} and \eqref{orbdafy} and indicate that  corresponding trajectories  move  clockwise on circles
\begin{equation}
\label{aaorb}
(x-1)^{2}+y^{2}=a^{2},\qquad a<1,
\end{equation}
of radius $a <1$ centred at $(1,0)$ when $c<1$ in \eqref{hprel}. 

Of special interest is the circle of unit radius for which $a=c=1$. 

Irrespective of prescribed initial conditions,  each member of the  family of orbits  created by  a sequence of different constants $c\geq 1$ possess the common property of lying in the positive half-plane and passing through the coordinate origin where they are mutually tangential. Dafermos \cite{d12} observes that global uniqueness  of solutions to \eqref{orb1}-\eqref{orb3} is violated on  allowing  switching at the origin between orbits of different entropy levels 
$H(x,y)$  defined by \eqref{hprel}.  Specifically,  a trajectory with initial conditions such that 
\begin{equation}
\label{dhinit}
H(x_{0},\,y_{0}):=\frac{(x^{2}_{0}+y^{2}_{0})}{2x_{0}}=c_{0}>1,
\end{equation}
on reaching the origin switches to another trajectory for which 
\begin{equation}
\label{dhsw}
H(x,y)=c_{1} >1,\qquad c_{1}\ne c_{0}.
\end{equation}

The composite trajectory is clearly continuous in the time variable $t$ and satisfies \eqref{dorb1} for almost all $t$; in fact for all $t$ for which $(x(t),\,y(t))\ne (0,\,0)$. From \eqref{orb2} the jump in the acceleration $\ddot{x}=\dot{y}$ at the origin is given by
\begin{equation}
\label{dj}
[\ddot{x}]_{(0,\,0)}=[\dot{y}]_{(0,\,0)}=c_{1}-c_{0}.
\end{equation}
Hence, $\dot{y}$ is bounded for all $t$ and $y$ is Lipschitz continuous.

Trajectories for which $c <1$ remain interior to the circle of radius $1$ and do not pass through the origin. Consequently, switching  is not possible.

The following theorem summarises for convenience these non-uniqueness results which are  also applicable to the general control problem \eqref{a1} and \eqref{a2}.  In  Section~\ref{ci}   the conclusion is compared to the non-uniqueness  theorem of  Buckmaster and Vicol \cite{bv19} for the Navier-Stokes equations.

\begin{figure}[htb]
\centerline{
\resizebox{6cm}{!}{
\def\firstcircle{(1,0) circle (1.0cm)}
\def\secondcircle{(1.5,0) circle (1.5cm)}
\begin{tikzpicture}[
  line cap=round,
  line join=round,
  >=Triangle,
  myaxis/.style={->,thick}]

\begin{scope}
        \begin{scope}[even odd rule]
            \clip \firstcircle (-3,-3) rectangle (3,3);
        \fill[pattern=north west lines,pattern color=blue] \secondcircle;
        \end{scope}
        \draw \firstcircle;
        \draw \secondcircle;
    \end{scope}

\node[above] at (3.2,0) {$R$};

\node[above] at (1.0,0) {$\Omega_2$};

\node at (2.2,1.8) {$\Omega_1$};

\draw [thick,-stealth](2.2,1.6) -- (2.2,0.8);

\draw[myaxis] (-0.66,0) -- (5.00,0) node[below left = 1mm] {$x$};      
\draw[myaxis] (0,-3.66) -- (0,3.66) node[below left = 1mm] {$y$};

\draw (1,0)--(1,-.1) node[below = 1mm] {$1$};
\end{tikzpicture}
}
}
\caption{\label{figomega}The region $\Omega=\Omega_{1}\cup\Omega_{2}$.}
\end{figure}
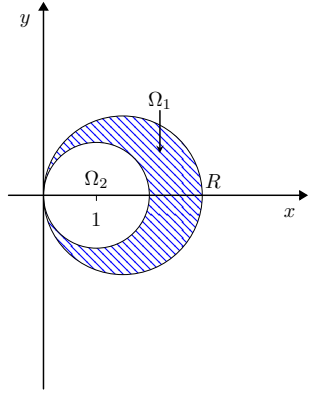

\begin{thm}
\label{dnonuniq}
For $R> 1$, define the plane region $\Omega =\Omega_{1}\cup\Omega_{2}$, depicted in Figure \ref{figomega},  by
\begin{eqnarray}
\label{O1def}
\Omega_{1}&:=& \left\{(x,y)\in\RR^{2}: 1\le (x-1)^{2}+y^{2},\, (x-R)^{2}+y^{2}\le R^{2}\right\},\\
\label{O2def}
\Omega_{2}&:=& \left\{ (x,y)\in\RR^{2}: (x-1)^{2}+y^{2}<1\right\}
\end{eqnarray}
and let $e(t)$ denote the piece-wise constant function
\begin{equation}
\label{eprof}
e(t)=\begin{cases}
                c_{0},& 0\le t\le t_{1},0\le t_{1}<2\pi,\\
                c_{n},& t_{n}=t_{1}+2(n-1)\pi\le t<t_{n}+2\pi,
\end{cases}
\end{equation}
for $n=1,2,\dots$ and $0\le t<\infty$.

Then for any sequence $c_{n}\geq 1$, there exists a Lipschitz continuous solution to \eqref{orb1} and \eqref{orb3}. Trivially, the initial value problem with initial data $(x(0)=x_{0},\,y(0)=y_{0})\in \Omega_{1} $ possesses an infinite number of Lipschitz continuous solutions. 
\end{thm}

Dafermos \cite{d12} regains uniqueness by suggesting  two admissibility criteria.  The first asserts that the physically meaningful solution to \eqref{dorb1} and \eqref{dorb2} on $\Omega=\Omega_{1}\cup\Omega_{2}$  is the limit as $\gamma\rightarrow 0^{+}$ of solutions $(x_{\gamma} ,\,y_{\gamma})$ to  the system with friction:
\begin{eqnarray}
\label{fric1}
\dot{x}_{\gamma}&=& y_{\gamma},\\
\label{fric2}
\dot{y}_{\gamma}&=& g(x_{\gamma},\,y_{\gamma}) -\gamma y_{\gamma},
\end{eqnarray}
where 
\[
g(x,y)=\begin{cases}
               \frac{(y^{2}-x^{2})}{2x},&(x-1)^{2}+y^{2}\ge 1,\\
                (1-x),& (x-1)^{2}+y^{2}<1.
               \end{cases}
\]

Solutions $(x_{\gamma},\,y_{\gamma})$ to \eqref{fric1} and \eqref{fric2} possess the limit $(x_{\gamma},\,y_{\gamma})\rightarrow (x,y)$ as $\gamma\rightarrow 0^{+}$ where
\begin{itemize}
\item[(i)]$(x_{\gamma}(0),\,y_{\gamma}(0))\rightarrow (x_{0},\,y_{0})$ and 
\begin{equation}
\nonumber
c_{0}=\frac{(x^{2}_{0}+y^{2}_{0})}{2x_{0}}.
\end{equation}
\item[(ii)] The limit $(x,\,y)$ moves clockwise on the circle radius $c_{0}$  until it reaches the origin $x=0,\,y=0$.
\item[(iii)] Upon reaching the origin, the solution switches \emph{once} to the circle with entropy  $e=1$ of radius $1$ centred at $(x,\,y)=(1,0)$. It remains on this circle moving clockwise.
\end{itemize}

\begin{figure}[htb]
\centerline{
\resizebox{6cm}{!}{
\begin{tikzpicture}[
  line cap=round,
  line join=round,
  >=Triangle,
  myaxis/.style={->,thick}]

\draw[myaxis] (0,0) -- (7.3,0) node[below left = 1mm] {$t$};      
\draw[myaxis] (0,0) -- (0,7) node[below left = 1mm] {$e(t)$};

\draw (0.3,0)--(0.3,-.1) node[below = 1mm] {$t_1$};
\draw (1.8,0)-- (1.8,-.1) node[below = 1mm] {$t_1+2\pi$};
\draw (3.3,0)--(3.3,-.1) node[below = 1mm] {$t_1+4\pi$};
\draw (4.8,0)--(4.8,-.1) node[below = 1mm] {$t_1+6\pi$};
\draw (6.3,0)-- (6.3,-.1) node[below = 1mm] {$t_1+8\pi$};

\draw (0,6)--(-.1,6) node[left] {$c_0$};
\draw (0,5)-- (-.1,5) node[left] {$c_1$};
\draw (0,3.5)--(-.1,3.5) node[left] {$c_3$};
\draw (0,2)--(-.1,2) node[left] {$c_2$};
\draw (0,1)-- (-.1,1) node[left] {$1$};

\draw[thick] (0,6)--(0.3,6);
\draw[thick] (0.3,5)--(1.8,5);
\draw[thick] (1.8,2)--(3.3,2);
\draw[thick] (3.3,1)--(4.8,1);
\draw[thick] (4.8,3.5)--(6.3,3.5);

\node at (0.2,6) {\tiny{$\times$}};
\node at (1.6,5) {\tiny{$\times$}};
\node at (3,2) {\tiny{$\times$}};
\node at (4.6,1) {\tiny{$\times$}};
\node at (6,3.5) {\tiny{$\times$}};

\draw (0.2,6)--(1.6,5);
\draw (1.6,5)--(3,2);
\draw (3,2)--(4.6,1);
\draw (4.6,1)--(6,3.5);

\end{tikzpicture}}
}
\caption{\label{disceprofile}Discrete energy profile.}
\end{figure}
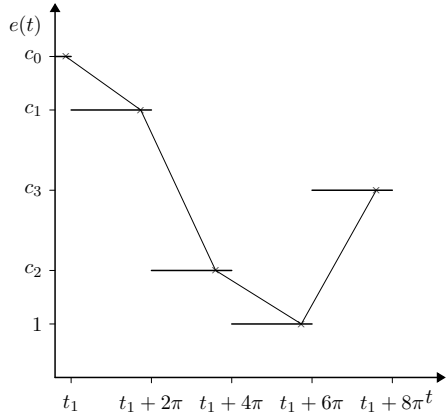

\begin{figure}[htb]
\centerline{
\resizebox{6cm}{!}{
\begin{tikzpicture}[
  line cap=round,
  line join=round,
  >=Triangle,
  myaxis/.style={->,thick}]

\draw[myaxis] (0,0) -- (7.3,0) node[below left = 1mm] {$t$};      
\draw[myaxis] (0,0) -- (0,7) node[below left = 1mm] {$e(t)$};

\draw (0,6)--(-.1,6) node[left] {$c_0$};

\draw (0,1)-- (-.1,1) node[left] {$1$};

\draw[thick] (0,6)--(1,6);
\draw[thick] (1,1)--(5,1);
\draw[thick,dashed] ((5,1)--(6.3,1);

\node at (0.8,6) {\tiny{$\times$}};
\node at (1.2,1) {\tiny{$\times$}};

\draw (0.8,6)--(1.2,1);

\end{tikzpicture}}
}
\caption{\label{disceprofilemax}Energy profile with maximal negative slope.}
\end{figure}
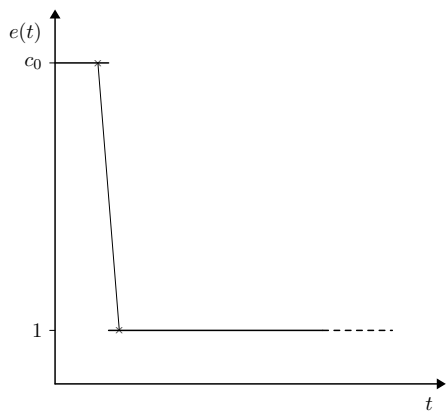

The second criterion proposed  by Dafermos is the \emph{entropy rate admissibility criterion}  which,  although not applied directly to the present problem, requires that physically admissible  solutions should have not only non-increasing entropy $e(t)$, but  that the entropy should decrease at maximum possible rate.  A simple illustration of the criterion  for the system \eqref{dorb1} and \eqref{dorb2} is  sketched in Figures \ref{disceprofile}-\ref{disceprofilemax} for  possible entropy profiles \eqref{eprof}. 

Let $\delta >0$ and in Figure \ref{disceprofile} select the point $\tau_{1}=t_{1}-\delta$ such that $0<\tau_{1}=t_{1}-\delta<t_{1}$ and $\tau_{1}$ lies within the first interval specified in \eqref{eprof} and therefore locates a point on the branch  $e=c_{0}$. 

 Subsequent points given by $\tau_{n}=t_{1}-\delta+2\pi n,\, n=1,2,\ldots$ are chosen to lie in the interval $(t_{n}, t_{n}+2\pi)$ so that $\tau_{n}$ is a point on $e=c_{n}$. On joining the points $\tau_{n},\,n=1,2,\ldots$ a piecewise linear graph is produced whose piecewise linear approximation is the original entropy function $e(t)$. The usual entropy criterion  stipulates decay of energy and accordingly the  piecewise linear graph must be  non-increasing. This, however, still permits an infinite number of profiles. A single preferred profile follows from the entropy rate admissibility criterion which identifies a single  profile corresponding to the  graph of  maximal negative slope  shown in Figure \ref{disceprofilemax}. Thus, the entropy rate criterion leads to the single solution derived  by the limiting friction argument. 

An immediate conclusion is that while it is possible for a smooth solution to \eqref{dorb1} and \eqref{dorb2} with  initial data in $\Omega_{1}$ to remain on the circle of radius $c_{0}$, the admissible unique solution given by either the limiting friction method or entropy rate criterion, though continuous, is  not smooth in the sense of \eqref{dj}. Consequently, the nonlinear oscillator considered by Dafermos contradicts the conjecture \emph{that to achieve uniqueness only smooth solutions should be admitted.} 

\begin{figure}[htb]
\centerline{
\begin{tikzpicture}[domain=0:3]
  \draw[->] (0,0) -- (3,0) node[right] {$t$};
  \draw[->] (0,0) -- (0,2) node[below left] {$e(t)$};
  \draw[thick,color=blue] plot (\x,{(1-tanh(3*(\x-1)))/2}); 
\end{tikzpicture}}
\caption{\label{smootheprofile}Smooth energy profile.}
\end{figure}
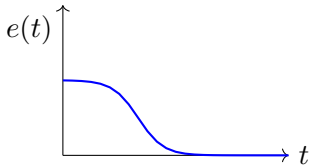

The notion of \emph{smoothness} is crucial not only here but for the  discussion in Section~\ref{ci}  of the entropy rate criterion in relation to the theorems to be there stated  of Buckmaster and Vicol \cite{bv19} and of De Lellis and Sz\`{e}kelyhidi \cite{ds09,ds10}. Consider, for example, the profile $e(t),\,e(0)\ne 0$ shown in Figure \ref{smootheprofile}. 

The entropy rate criterion selects the preferred profile as that having maximal decreasing derivative $e'(t)$. Consequently, for any fixed energy profile $e(t)$ there will always be one with a more rapidly decreasing profile. Hence all non-zero smooth profiles given by the results of Buckmaster and Vicol \cite{bv19} and by De Lellis and  Sz\`{e}kelyhidi \cite{ds10,ds13} are inadmissible according to the entropy rate criterion. But clearly the limit of such profiles is the maximal profile $e(t)=e(0)H(t)$ where $H(t)$  satisfies 
 $H(0)=1$ and $H(t)=0$ for $t>0$. This maximal profile is neither  smooth nor even a weak solution. 
A similar argument  used by Feireisl \cite{ef14}  proves inadmissibility of  weak solutions to the compressible Euler equations derived by convex integration methods.  On the other hand, when $e(0)=0$, the maximal entropy rate  yields the trivial admissible profile $e(t)=0,\,t\ge 0$. This remark indicates that the profile found by Scheffer \cite{s93} and by Shnirelmann \cite{s97} for the Euler equations is  inadmissible. 

\begin{rem}[Regular oscillations]
As noted by Dafermos \cite{d12}, the  entropy rate admissibility criterion applied to the nonlinear oscillator implies  that starting from any initial data, the oscillation  after once having passed through the origin  subsequently  is  regularly  periodic for all time.
\end{rem}

The next task is to interpret the solution to \eqref{orb1} and \eqref{orb3} as Lagrange trajectories appropriate to certain  steady plane fluid flows satisfying either the Euler or Navier-Stokes equations.


\section{Plane Euler and Navier-Stokes equations}\label{ens}
This section  examines implications of  Theorem~\ref{dnonuniq} for compressible steady flows governed by  Euler and Navier-Stokes equations on the plane region $\Omega_{1}$ {exterior to the unit circle defined by \eqref{O1def} and on the region $\Omega=\Omega_{1}\cup\Omega_{2}$}.  (See Figure \ref{figomega}.) It is supposed that \eqref{dorb1} represents  Lagrange trajectories of fluid particles and that for both the Euler and Navier-Stokes  equations the velocity vector has cartesian components
\begin{eqnarray}
\label{eu1}
u(x,y)&=& y,\\
\label{eu2}
v(x,y)&=& \frac{(y^{2}-x^{2})}{2x},
\end{eqnarray}
 for which the entropy $H(x,y)$ defined by \eqref{hprel}  is constant.

{\begin{rem}[Lagrange trajectories and Euler's equation.] \label{remark31}
Throughout this Remark, vector quantities are denoted  in bold type. Lagrange trajectories and the Euler equations represent two distinct methods for describing fluid motion. Lagrange trajectories trace the motion of a fluid particle in space and time and consequently the vector velocity field $\textbf{v}$ is a function of initial spatial position $\textbf{x}_{0}$ and time $t$. The position vector $\textbf{x}(t)=\textbf{x}(\textbf{x}_{0},t)$ and velocitiy vector $\textbf{v}(\textbf{x},t)=\textbf{v}(\textbf{x}_{0},t)$ are related by the  system of ordinary differential equations:
\begin{equation}
\label{laode}
\dot{\textbf{x}}=\textbf{v}(\textbf{x}_{0},t),\qquad \textbf{x}(0)=\textbf{x}_{0},
\end{equation}
 and for a solution to exist and to  be unique it is sufficient that $\textbf{v}$ is  Lipschitz continuous. There are, however, conditions under which the solution may neither exist nor be unique.

The same velocity field $\textbf{v}$ occurs in the Euler description except that  the motion is considered  at fixed spatial position  $\textbf{x}$ and at  varying time. The governing equations are the Euler system of partial differential equations subject to prescribed initial and boundary data. Existence and uniqueness of a solution follow from theorems in partial differential equations and depend upon function spaces and definition of solution.  The solution, however defined, may be non-unique, 
but can  be used as the velocity in the ordinary differential equations \eqref{laode} to derive the Lagrange trajectories which likewise may be non-unique   for given initial data. The procedure may be reversed and the possible non-unique  Lagrange trajectories first computed for specified velocity field. This velocity is then substituted in the Euler equations and the corresponding pressure, density, and initial and boundary data  calculated.  The initial boundary problem obtained by this semi-inverse method may be non-unique  and    solutions may exist additional to the velocity used in the system \eqref{laode}.

The interrelation between solutions for Lagrange trajectories and for the Euler system is remarked upon by Robinson and Sadowky \cite{rs09a, rs09b}.
The semi-inverse procedure is  adopted  in what follows.


Similar remarks apply to the connexion between Lagrange trajectories and the Navier-Stokes equations  examined later in this Section.

It is noted that  Lagrange trajectories have recently found applications in oceanography, atmospherics and biology. \end{rem}
}

The following lemmas establish that the particular component velocities on the right of \eqref{eu1} and \eqref{eu2}  satisfy the compressible Euler equations with  mass density
\begin{equation}
\label{edens}
  \rho(x)=x^{-1}.
\end{equation}

 {
\begin{rem}
The conclusion, however, is not confined to  velocity components \eqref{eu1} and \eqref{eu2}, nor to the mass density \eqref{edens}. Appendix~\ref{A3} explains how a certain   general class of velocities derived from a conservative system of Lagrange trajectories also satisfies  Euler's equations for suitable choice of mass density and pressure.
\end{rem}
}

\begin{lem}\label{econt}
In the interior of $\Omega_{1}$, defined by \eqref{O1def}, the velocity components \eqref{eu1} and \eqref{eu2} satisfy the continuity equation
\begin{equation}
\label{ceqn}
\frac{\partial}{\partial x}\left(\rho u\right)+\frac{\partial}{\partial y}\left(\rho v\right)=0.
\end{equation}
\end{lem}
\noindent \textbf{Proof}: By direct substitution. Note, however, that since
\begin{equation}
\label{eucomp}
\frac{\partial u}{\partial x}+\frac{\partial v}{\partial y}= \frac{y}{x}, \qquad x\neq 0,
\end{equation}
the fluid is compressible except possibly on $y=0,\, x\ne 0.$


\begin{lem}\label{ebc}
On $\partial\Omega_{1}\backslash \left\{(0,0)\right\}$, the boundary condition
\begin{equation}
\label{ebc}
 (un_{1}+ vn_{2})=0,\qquad (x,y)\in \partial\Omega_{1}\backslash (0,0),
\end{equation}
holds, where $n_{1},\,n_{2}$ are the cartesian components of any normal vector on the boundary. 
\end{lem}

\noindent \textbf{Proof}:  
The inner boundary  given by
\begin{equation}
\label{inbdr}
(x-1)^{2}+y^{2}=1, \qquad (x,y)\neq (0,0),
\end{equation}
has  unit normal whose components are   $\left((x-1),\,y\right)$. Moreover, on \eqref{inbdr} 
\begin{equation}
\nonumber
u(x,y)=y,\qquad v(x,y) =(1-x),
\end{equation}
and  \eqref{ebc} is immediate when $(x,y)\neq (0,0)$. The same argument applies to the outer boundary
\begin{equation}
\label{outer}
(x-R)^{2}+y^{2}=R^{2},\qquad R>1.
\end{equation}

Condition \eqref{ebc} trivially implies   tangential flow at the boundary.

\begin{lem}\label{emtm}
The given velocity and density on the interior of $\Omega_{1}$ satisfy the  balance of steady  linear momentum 
\begin{eqnarray}
\label{elin1}
\frac{\partial}{\partial x}\left(\rho u^{2}\right)+\frac{\partial}{\partial y}\left(\rho u v\right)+\frac{\partial p}{\partial x}&=& 0,\\
\label{elin2}
\frac{\partial}{\partial x}\left(\rho u v\right)+\frac{\partial}{\partial y}\left(\rho v^{2}\right)+\frac{\partial p}{\partial y}&=&0,
\end{eqnarray}
subject to the specified pressure
\begin{equation}
\label{epress}
p(x,y):= \frac{(x^{2}+y^{2})}{2x}-1  = H-1.
\end{equation}
\end{lem}
\noindent As $H$ is constant along trajectories, we hence have the pressure $p$ constant along trajectories.\\

\noindent \textbf{Proof}: Insertion of velocity and density into the balance of steady linear momentum \eqref{elin1}   after integration  yields
\begin{equation}
\nonumber
p(x,y)= \frac{(x^{2}+y^{2})}{2x}+f(y),
\end{equation}
where $f$ is an arbitrary function.  On the other hand, \eqref{elin2} leads to 
\begin{equation}
\nonumber
p(x,y)= \frac{y^{2}}{2x}+g(x),
\end{equation}
for arbitrary function $g$.
Set $f\equiv 0$ and $g=(x/2+ b)$ for arbitrary constant $b$ taken to be $b=-1$ to ensure that $p(x,y)$ vanishes on the inner boundary \eqref{inbdr}. 
\vspace{0.2cm}

Lemma~\ref{econt}, Lemma~\ref{ebc} and Lemma~\ref{emtm} establish the following theorem,

\begin{thm}\label{ethm}
When the velocity, density and pressure satisfy \eqref{eu1}, \eqref{eu2}, \eqref{edens}, and \eqref{epress}, the fluid flow with Lagrange particle trajectories \eqref{dorb1} and \eqref{dorb2} in $\Omega_{1}$ satisfy the continuity equation \eqref{ceqn},  balance of steady linear momentum \eqref{elin1} and \eqref{elin2}, and  the boundary condition \eqref{ebc}.
\end{thm}

  Note that the velocity field \eqref{eu1} and \eqref{eu2} generates  Lagrange trajectories that are non-unique for  initial data specified  in $\Omega_{1}$.  A unique Lagrange trajectory is obtained from the entropy rate admissibility criterion employed in accordance with the construction of the previous section.  

The region $\Omega_{1}$ occupied by the fluid may be enlarged  to include the interior of the unit circle $(x-1)^{2}+y^{2}=1$ as proposed by Dafermos; see  \eqref{dorb2}. As before (see \eqref{O2def}), denote the unit disc by
\begin{equation}
\label{om2def}
\Omega_{2}:=\left\{(x,y): (x-1)^{2}+y^{2}\le 1\right\}.
\end{equation}
Then  for $(x,y)\in \Omega_{2}$  put 
\begin{eqnarray}
\label{rhodisc}
\rho(x,y)&=& 1,\\
\label{udisc}
u(x,y)&=& y, \\
\label{vdisc}
v&=& (1-x).
\end{eqnarray}

It follows from  mass conservation, or directly, that 
\begin{equation}
\label{incdisc}
\left(\frac{\partial u}{\partial x}+\frac{\partial v}{\partial y}\right)=0,\qquad (x,y)\in \Omega_{2},
\end{equation}
and the flow is incompressible.  Expressions for balance of  steady linear momentum become
\begin{eqnarray}
\nonumber
\frac{\partial}{\partial x}(u^{2})+\frac{\partial}{\partial y}(uv)+\frac{\partial p}{\partial x}&=&0, \qquad (x,y)\in \Omega_{2},\\
\nonumber
\frac{\partial}{\partial x}(uv)+\frac{\partial}{\partial y}(v^{2})+\frac{\partial p}{\partial y}&=&0, \qquad (x,y)\in \Omega_{2},
\end{eqnarray}
and  are satisfied for   pressure $p(x,y)$  given by 
\begin{equation}
\label{presdisc}
p(x,y)= \frac{(x^{2}+y^{2}-2x)}{2}, \qquad (x,y)\in \Omega_{2},
\end{equation}
which vanishes on $\partial\Omega_{2}$. In consequence, the pressure in $\Omega=\Omega_{1}\cup\Omega_{2}$ is continuous across the inner  boundary $\partial\Omega_{2}$ (the  unit circle \eqref{inbdr}) except at the origin $(0,\,0)$. In this respect, an arbitrary constant can be added to the pressure \eqref{epress} inside $\Omega_{1}$ provided the same constant is added to the pressure \eqref{presdisc} in $\Omega_{2}$.

 The vector field \eqref{udisc} and \eqref{vdisc} on the unit circle satisfies the boundary condition
\begin{equation}
\label{bcdisc}
(un_{1}+vn_{2})=0,\qquad (x,y)\in \partial\Omega_{2},
\end{equation}
and fluid particles inside $\Omega_{2}$ flow tangential to $\partial \Omega_{2}$ and cannot penetrate into  $\Omega_{1}$. Equally, particles flowing inside $\Omega_{1}$ can never reach inside $\Omega_{2}$, but as already shown, flow tangential to $\partial \Omega_{2}$. It may be easily checked by direct computation that the Rankine-Hugoniot jump conditions are satisfied across the inner boundary $\partial\Omega_{2}.$

The Lagrange trajectories in $\Omega_{2}$  are circles centred at $(1,0)$ and of radius $c,\,0\le c \le 1$ and are uniquely defined by initial data. An admissibility criterion is therefore required only for trajectories in the region $\Omega_{1}$ and as previously shown leads to the unique trajectory in which particles steadily traverse the unit cirle \eqref{inbdr}.  Fluid behaviour in the composite region $\Omega=\Omega_{1}\cup\Omega_{2}$ subject respectively to the velocities \eqref{eu1},\eqref{eu2}, \eqref{udisc}, \eqref{vdisc} and pressures \eqref{epress} and \eqref{presdisc} in $\Omega_{1}$ and $\Omega_{2}$ results in the fluid in $\Omega_{1}$ rapidly becoming  a steady swirling motion entirely  confined clockwise to the unit circle \eqref{inbdr}. The inner region $\Omega_{2}$ remains fully occupied by incompressible fluid moving with velocity \eqref{udisc} and \eqref{vdisc} subject to pressure \eqref{presdisc}.

\vspace{0.5cm}
Similar conclusions to Theorem~\ref{ethm} are valid for the plane compressible steady flow satisfying the  Navier-Stokes equations.

\begin{thm}\label{nsthm}
Let the velocity, density and pressure satisfy \eqref{eu1}, \eqref{eu2}, \eqref{edens}, and \eqref{epress} in the region $\Omega_{1}$ .  Then for  viscosity coefficients $\mu,\,\lambda$ 
given by
\begin{eqnarray}
\label{mu}
\mu&=& 4\cos^{2}{\theta},\\
\label{lamb}
\lambda&=& -4\theta \cot{\theta}-4\cos^{2}{\theta},
\end{eqnarray}
where $x=r\cos{\theta},\,y=r\sin{\theta}$, the particle trajectories \eqref{dorb1} in $\Omega_{1}$ are solutions to the Navier-Stokes equations specified by the continuity equation \eqref{ceqn}, the balance of steady linear momentum
\begin{eqnarray}
\label{ns11}
\frac{\partial(\rho u^{2})}{\partial x}+\frac{\partial(\rho uv)}{\partial y}&=&
\frac{\partial\Sigma_{11}}{\partial x}+\frac{\partial\Sigma_{21}}{\partial y},\\
\label{ns22}
\frac{\partial (\rho uv)}{\partial x}+\frac{\partial (\rho v^{2})}{\partial y}&=& 
\frac{\partial \Sigma_{21}}{\partial x}+\frac{\partial \Sigma_{22}}{\partial y},
\end{eqnarray}
where
\begin{eqnarray}
\label{ns2}
\Sigma_{\alpha\beta}&=& -p\delta_{\alpha\beta}+\sigma_{\alpha\beta},\qquad \alpha,\,\beta =1,2,\\
\label{ns11}
\sigma_{11}&=& (\lambda +2\mu)\frac{\partial u}{\partial x}+\lambda \frac{\partial v}{\partial y},\\
\label{ns22}
\sigma_{22}&=& (\lambda +2\mu)\frac{\partial v}{\partial y}+\lambda \frac{\partial u}{\partial x},\\
\label{ns12}
\sigma_{12}&=& \mu\left(\frac{\partial u}{\partial y}+\frac{\partial v}{\partial x}\right),
\end{eqnarray}
and the boundary condition \eqref{ebc}. Here, $\delta_{\alpha\beta}$ denotes  the standard Kronecker delta.
\end{thm}

\noindent \textbf{Proof}: By direct substitution of the stated expressions.

\vspace{0.3cm}

Theorem~\ref{nsthm} has similar implications to Theorem~\ref{ethm}  for the relationship between Lagrange trajectories and the Navier-Stokes equations. A given velocity that solves the compressible Navier-Stokes equations in $\Omega_{1}$ leads to  corresponding non-unique Lagrange trajectories  for the same initial data. The entropy rate admissibility criterion, however, distinguishes  a unique trajectory.

By comparison, it is shown by Robertson and Sadowkski \cite{rs09a,rs09b}  for the incompressible Navier-Stokes equations on a bounded three-dimensional region subject to Dirichlet boundary data, that for sufficiently regular velocities the Lagrange trajectories are unique for almost all initial data.

As before, the analysis may be extended to the enlarged region $\Omega=\Omega_{1}\cup\Omega_{2}$ employing \eqref{udisc} and \eqref{vdisc} as velocity components in the unit disc $\Omega_{2}$. Besides the incompressibility condition \eqref{incdisc}, the gradient of the velocity components \eqref{udisc} and \eqref{vdisc}   for $(x,y)\in\Omega_{2}$ satisfy the relations
\begin{equation}
\label{nordisc}
\frac{\partial u}{\partial x}=\frac{\partial v}{\partial y}= 0,
\end{equation}
and therefore $\sigma_{\alpha\beta}=0$ for $\alpha,\, \beta=1,2.$ That is, the viscous contribution vanishes identically, and the Navier-Stokes equations are satisfied for all viscosity coefficients $\lambda,\,\mu$ inside the unit disc $\Omega_{2}$. Accordingly, Lagrange trajectories in $\Omega_{2}$ exhibit the same properties as those noted for the Euler equations implying that the fluid particles  in $\Omega_{1}$ and $\Omega_{2}$ cannot interpenetrate. The entropy rate admissibility criterion  establishes that the unique motion in $\Omega_{1}$ is again concentrated solely on the unit circle around which it continuously swirls clockwise.

 A final important observation is that in the problems of this section  the density $\rho(x)=x^{-1}$ in $\Omega_{1}$ becomes singular on the $y-$axis leading to density blow-up. Fluid particle aggregation is therefore to be expected and  occurs for the entropy rate admissible trajectories  but not for the ``wild'' non-admissible ones. In fact, this expectation is reflected in numerical computations presented in the next section.

\section{Numerical approximation}\label{numapp}
This section discusses the numerical approximation of solutions to the system \eqref{dorb1} and \eqref{dorb2} relevant for Lagrange trajectories in the combined region $\Omega=\Omega_{1}\cup\Omega_{2}$.  Of special interest are numerical methods which preserve the qualitative behaviour of the flow in phase space and  induce convergence to the solution selected by the entropy rate admissibility criterion. The basic example of such geometric numerical integrators (see \cite{hlw06}) is the symplectic Euler method applied here.
 A constant step size $h>0$ is used to compute approximations $(x_{n},\,y_{n})$ of the solution $(x(nh),\,y(nh)),\, n\in N   $ to \eqref{dorb1} and \eqref{dorb2}. As before, denote initial conditions by $(x_{0},\,y_{0}):=(x(0),\,y(0))$  and define $x_{n},\,y_{n}$ by 
\begin{equation}\begin{rcases}  y_{n+1} = y_n+h\frac{y_n^2-x_n^2}{2 x_n},\qquad x_{n+1} = x_n + hy_{n+1},  &\text{for } (x_n-1)^2+y_n^2\geq 1,\\ 
  y_{n+1} = y_n+ h(1-x_n),   \quad     x_{n+1} = x_n+h y_{n+1}, & \text{for } (x_n-1)^2+y_n^2< 1.
\end{rcases}\label{num}\end{equation}
As described in \cite[Chapt.~IX]{hlw06}, the symplectic Euler method determines a flow that preserves an approximation to the entropy function $H(x,y)$, defined in \eqref{hpprel},  to an error of order $h$.

\begin{figure}[htb]
\centerline{\resizebox{10cm}{!}{\includegraphics{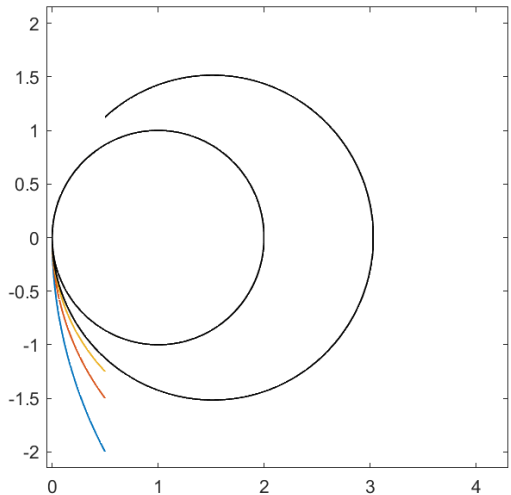}}}
\caption{\label{figinitial}Trajectories of symplectic Euler discretisation for different initial conditions.}
\end{figure}

Figure \ref{figinitial} illustrates trajectories of \eqref{num} corresponding to  four initial conditions $(x_{0},\,y_{0})$ satisfying $(x_{0}-1)^{2}+y_{0}^{2}>1$, when $h=10^{-4}$. In all cases, the numerical trajectories traverse the level sets of $H$ to high accuracy. Near the origin $(0,\,0)$, the trajectories switch to the unit circle $(x-1)^{2}+y^{2}=1$ as expected from the entropy rate admissibility criterion. The black trajectory, which corresponds to the initial condition $(x_{0},\,y_{0})=(\frac{1}{2},\,\frac{9}{8})$, is further analysed  below.

\begin{figure}[htb]
\centerline{\resizebox{7cm}{!}{\includegraphics{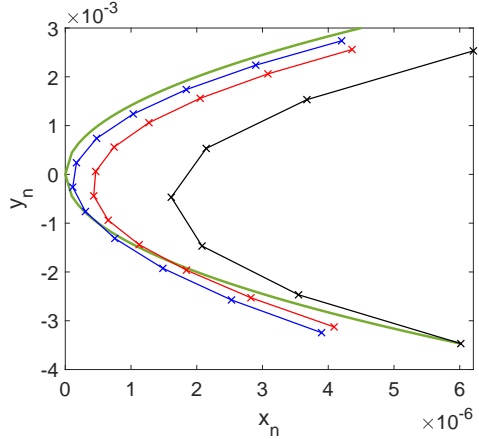}}}
\caption{\label{figtraj}Trajectories for symplectic Euler discretisation near $(0,\,0)$ for $h=2.5\times 10^{-4}$(blue); $h=5\times 10^{-4}$(red); $h=10^{-3}$(black). The unit circle $(x-1)^{2}+y^{2}=1$ in green.}
\end{figure}

Figure \ref{figtraj} shows the behaviour near the origin in varying detail depending upon step size: $h=2.5\times 10^{-4}$(blue),  $h=5\times 10^{-4}$(red),  $h =10^{-3}$(black). All trajectories start from the same initial condition $(x_{0},\,y_{0})$,  cross the $x-$axis  at a positive value $x\sim h^{2}$ and then follow the unique solution close to the unit circle. The well-known analysis of the symplectic Euler method of the linear  equation in the domain  $(x-1)+y^{2} <1$ shows that the trajectory remains in a neighbourhood of the unit circle of order $h$. The rigorous analysis presented in \cite[Chapt.~11]{hlw06}, however, is required to establish that trajectories always cross the $x-$axis at points $x >0,$ and thereafter enter and remain near the unit circle.  In fact, for small $h$, the level sets of $H$ are perturbed to the right near the origin. The geometry of the discrete flow explains the selection  of the solution obtained  by the entropy rate admissibility criterion.

\begin{figure}[htb]
\centerline{\resizebox{7cm}{!}{\includegraphics{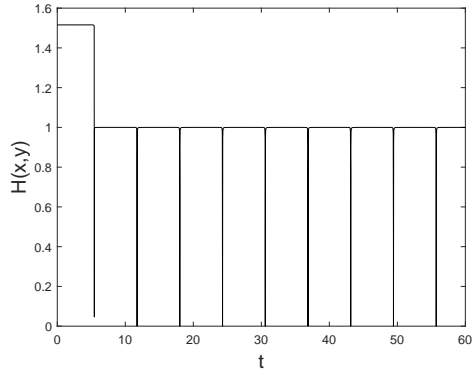}}}
\caption{\label{figenergy} $H(x,y)$ as function of $t$ for initial condition $(x_0,\,y_0) = (\frac{1}{2},\frac{9}{8})$.}
\end{figure}

\begin{figure}[htb]
\centerline{\resizebox{7cm}{!}{\includegraphics{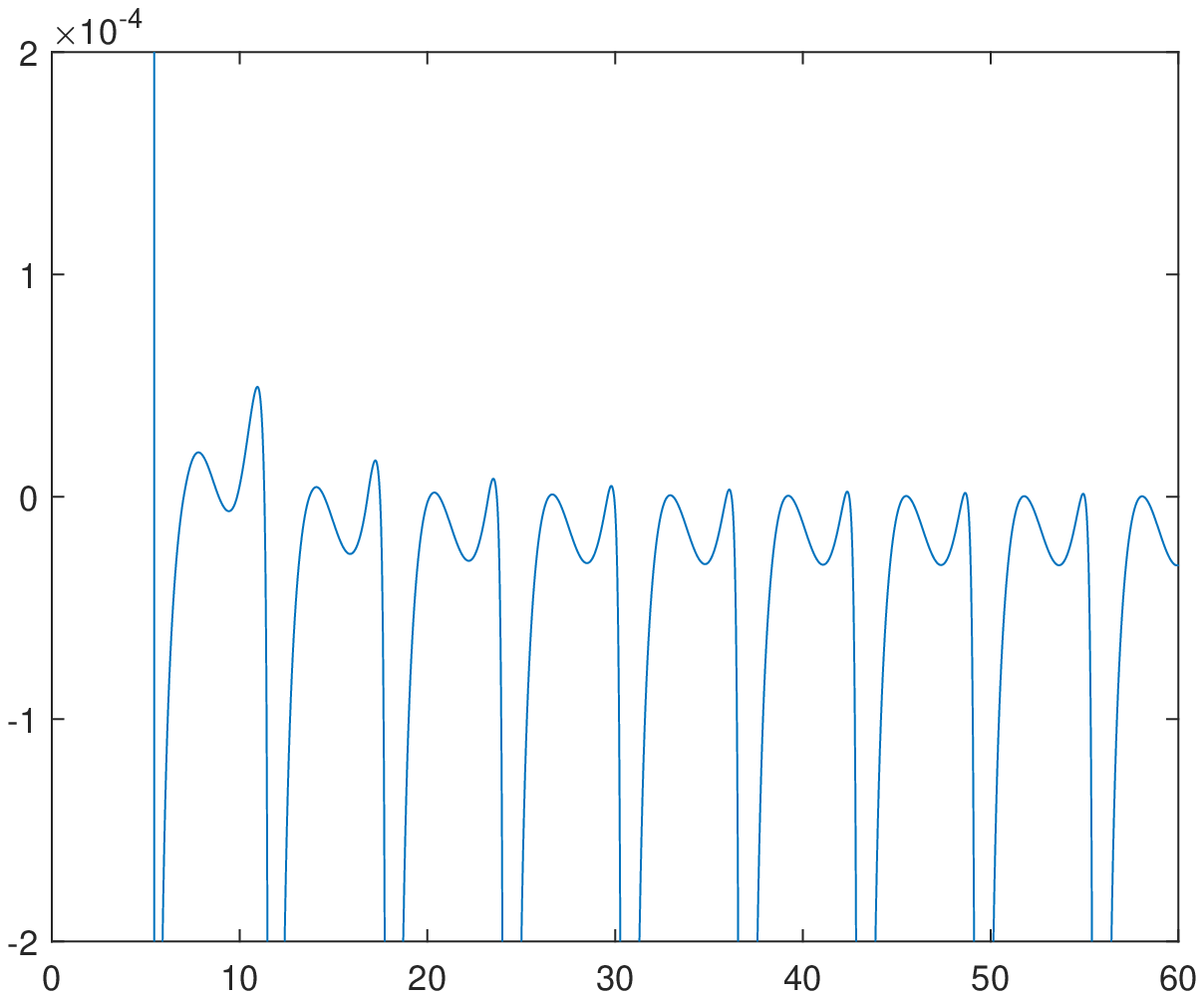}}}
\caption{\label{figerror}Error $(H(x,y)-1)$ as function of $t$ for initial condition $(x_{0},\,y_{0})= (\frac{1}{2},\,\frac{9}{8}).$}
\end{figure}

\begin{figure}[htb]
\centerline{\resizebox{7cm}{!}{\includegraphics{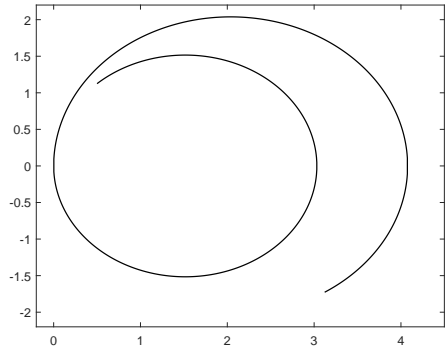}}}
\caption{\label{figexpeul}Trajectory of explicit Euler discretisation for initial condition $(x_{0},\,y_{0})= (\frac{1}{2},\,\frac{9}{8}).$}
\end{figure}

Figure \ref{figenergy} depicts the numerically computed entropy profile $t\mapsto H(x(t),\,y(t))$ for initial condition $(x_{0},\,y_{0})=(\frac{1}{2},\,\frac{9}{8})$ introduced above. Again take  $h=10^{-4}$. When the trajectory first passes through the origin, $H$ drops to a  value near $1$ and then remains approximately constant.

Figure \ref{figerror} indicates the discretisation error in the entropy  profile $t\mapsto \left[H(x(t),\,y(t))-1\right].$  The error is of size less than $10^{-4}= O(h)$ except for floating point errors when $(x,\,y)$ is near the origin where  $H$  possesses a singularity.

It has  thus been shown how geometric properties in phase space of the symplectic Euler method select the solution to the system \eqref{dorb1} and \eqref{dorb2} according to the entropy rate admissibility criterion. Numerical methods without such properties are not expected to select these solutions. 

Figure \ref{figexpeul} shows the trajectory computed by the explicit Euler method again for  initial condition  $(x_{0},\,y_{0})=(\frac{1}{2},\,\frac{9}{8})$. For this discretisation, the numerical solution moves to a larger circle with increased entropy $H$, even with smaller step size $h=10^{-5}.$

\section{Relevant convex integration results}\label{ci}

Section~\ref{mot} describes how the entropy rate admissibility criterion applied to the constant functions \eqref{hprel} selects a unique orbit from the infinitely many circular paths generated by  the Dafermos equations \eqref{orb1} and \eqref{orb3}. A piecewise non-increasing linear graph of maximal negative slope leads to the requisite single profile. A notable feature of the construction, stated in Theorem~\ref{dnonuniq}, is the uncountable number of profiles that can be arbitrarily chosen for prescribed initial conditions.

This feature resembles results of Buckmaster and Vicol \cite{bv19} and of De Lellis and Sz\`{e}kelyhidi \cite{ ds13} which for convenience are recalled in the following theorem.

\begin{thm}\label{bv}
There exists $\beta >0$ such that for any non-negative smooth function $E(t):[0,T]\rightarrow \RR_{\ge 0}$ there exists a weak vector solution $\textbf{v}\in C^{0}_{t}\left([0,T];H^{\beta}_{x}(\TT^{3})\right)$ of the Navier-Stokes equations such that 
\begin{equation}
\nonumber
\int_{\TT^{3}}|\textbf{v}(x,t)|^{2}\,dx=E(t),
\end{equation}
for all $t\in [0,T]$. Moreover, the associated vorticity $\nabla\times \textbf{v}$ lies in $C^{0}_{t}\left([0,T];L^{1}_{x}(\TT^{3})\right),$
where $\TT^{3}=\RR^{3}\backslash 2\pi Z^{3}$ denotes the unit cube with periodic boundary conditions.
\end{thm}

It is apparent that arbitrary prescription of the  energy $E(t)$  produces  an  infinite number of solutions  in the given class for the same initial conditions.  An admissibility criterion is required  to identify a single physically relevant velocity.  Whether an entropy rate admissibility criterion of the type considered here is appropriate in this respect remains an important open question.

A detailed comparison with the analysis of Section~\ref{ens} reveals distinct differences and any analogy between the two treatments is likely to be superficial.  Indeed,  the entropy  $H(x,\,y)$ used in Section~\ref{ens} is not a local (kinetic) energy function $\eta(t)$ defined by
\begin{equation}
\label{edef}
\eta(t):=\frac{\rho}{2} \left(\dot{x}^{2}+\dot{y}^{2}\right)
\end{equation}
where $\rho(x,y)$ is a  specified density. The relation to $H(x,y)$ is obtained by substitution of \eqref{a1} and \eqref{a2} in \eqref{edef} and is given by
\begin{equation}
\label{eH}
\eta(t) = 2x^{2}w^{2}(x,y)\rho(x,y) H^{2}(x,y).
\end{equation}
Consequently,  functions $\eta$ and $H^{2}$ are identical only when $2x^{2} w^{2}\rho=1$.


Finally, note that   uniqueness properties arising in convex integration studies refer to  solutions of the respective equations and not to the associated Lagrange trajectories.  Whether there are non-unique trajectories for certain  non-unique convex integration solutions remains a second open question.


\paragraph{Acknowledgements:}\emph{ Valuable discussions during the preparation of this paper  with Professors Z. Artstein, E. Feireisl, A. Ostermann  and E. Titi are gratefully acknowledged. Topics treated in this paper first arose during the 2021 Workshop on \emph{Convex Integration and Partial Differential Equations} organised by ICMS, Edinburgh.}

\bibliographystyle{amsplain}

\section*{Appendices}\label{app}

\appendix

\section{State variables.}\label{orbcal}

 Let state variables $(x(t),\,y(t))$ satisfy \eqref{orb1} and \eqref{orb3} in region $\Omega_{1}$ defined by \eqref{O1def}.  Consider for $x\ne 0,\,y\ne 0,$ the expression 
\begin{eqnarray}
\nonumber
\frac{d}{dt}\left(\frac{y^{2}}{x}\right)&=& \frac{y}{x^{2}}\left(2\dot{y}x-y\dot{x}\right)\\
\nonumber
&=& \frac{y}{x^{2}}\left(y^{2}-x^{2} -y\dot{x}\right)\\
\nonumber
&=& \frac{y}{x^{2}}\left(y\dot{x}-x^{2}-y\dot{x}\right)\\
\nonumber
&=& -y\\
\nonumber
&=& -\dot{x}, 
\end{eqnarray}
which on integration gives (cp., \eqref{codef}) 
\begin{equation}
\label{aic}
\frac{y^{2}}{x}+x= \frac{y^{2}_{0}}{x_{0}}+x_{0}\equiv 2c_{0}
\end{equation}
or
\begin{equation}
\nonumber
(x-c_{0})^{2}+y^{2}= c_{0}^{2},
\end{equation}
which may be rewritten as
\begin{equation}
\nonumber
y(t) = \pm\sqrt{x(2d-x)}\qquad  (= \dot{x}).
\end{equation}
On setting $c_{0}z(t)=(x(t)-c_{0})$ the last expression becomes
\begin{equation}
\nonumber
\dot{z}(t)=\pm \sqrt{(1-z^{2})},
\end{equation}
so that 
\begin{equation}
\nonumber
z(t)=\pm\cos{(\mp t+\sin^{-1}{z_{0}})},
\end{equation}
and therefore 
\begin{equation}
\label{ax}
x(t)=c_{0}\pm c_{0}\cos{(\mp t+\sin^{-1}{z_{0}})}.
\end{equation}
  A corresponding expression for $y(t)$ obtained on differentiation of the last relation is
\begin{equation}
\label{ay}
y(t)=\pm c_{0}\sin{(\mp t+\sin^{-1}{z_{0}})}.
\end{equation}
The positive root in \eqref{ax} and \eqref{ay} yields \eqref{orbdafx} and \eqref{orbdafy} with $\theta_{0}=\sin^{-1}{z_{0}}$ given by (cp., \eqref{thdef}) 
\begin{equation}
\nonumber
\tan{\theta_{0}} =\frac{y_{0}}{(x_{0}-c_{0})}.
\end{equation}

Alternatively, elimination of the dependent variable $y$   between \eqref{orb1} and \eqref{orb2} leads to a differential equation for $x$ whose solution is 
\begin{equation}
\label{dafx}
x(t)-c_{0}= A\cos{(t)}+B\sin{(t)},
\end{equation}
and consequently
\begin{equation}
\label{dafy}
y(t)=-A\sin{(t)}+B\cos{(t)},
\end{equation}
where $A,\,B,$ constants determined by initial conditions $(x_{0},\,y_{0})$,   are given by
\begin{eqnarray}
\label{Adef}
A&=& x_{0}-c_{0},\\
\label{Bdef}
B&=& y_{0},
\end{eqnarray}
and $c_{0}$ is specified by \eqref{aic}. Expressions \eqref{dafx} and \eqref{dafy} are equivalent to \eqref{ax} and \eqref{ay}.



\section{Navier-Stokes equations.}\label{nsapp1}
Throughout  Appendix~\ref{nsapp1} and Appendix~\ref{A2}  a suffix notation and summation convention are adopted together with a subscript comma to denote partial differentiation. Components of vectors, unless otherwise stated, are with respect to a given Cartesian coordinate system.

Consider a compressible fluid (of variable density) in plane steady motion that satisfies the isotropic Navier-Stokes equation for a given velocity field having  Cartesian components $(v_{1},\,v_{2})$. The aim, using a semi-inverse method, is to determine the viscosity  coefficients $\lambda(x_{1},\,x_{2})$ and $\mu(x_{1},\,x_{2})$ that ensure the corresponding stress distribution is in equilibrium subject to zero body force. The region occupied by the fluid is not defined at this stage.  Nor are boundary conditions considered.

The compatible strain rate  components, given by
\begin{equation}
\label{strai}
e_{\alpha\beta}=\frac{1}{2}\left(v_{\alpha,\beta}+v_{\beta,\alpha}\right),\qquad \alpha,\,\beta =1,2,
\end{equation}
are related to the stress components by \eqref{ns11}-\eqref{ns12}  concisely 
 written as
\begin{equation}
\label{ssrel}
\sigma_{\alpha\beta}=\lambda e_{\gamma\gamma}\delta_{\alpha\beta}+2\mu e_{\alpha\beta},
\end{equation}
where  $\delta_{\alpha\beta}$ denotes the Kronecker delta. The viscosity  coefficients $\lambda(x_{1},\,x_{2}),\,\mu(x_{1},\,x_{2})$ are functions 
of position to be determined such that the stress is in equilibrium under zero body force; that is 
\begin{equation}
\label{sequi}
\sigma_{\alpha\beta,\beta}=0.
\end{equation}
In consequence, the Navier-Stokes equations \eqref{ns11} and \eqref{ns22} reduce to the Euler equations \eqref{elin1} and \eqref{elin2}.

Easy deductions from \eqref{ssrel} are the trace relation
\begin{equation}
\label{trace}
\sigma_{\alpha\alpha}=2(\lambda+\mu)e_{\alpha\alpha},
\end{equation}
and
\begin{equation}
\label{fundrel}
2\mu=\frac{\left(\sigma_{11}-\sigma_{22}\right)}{\left(e_{11}-e_{22}\right)}=\frac{\sigma_{12}}{e_{12}},
\end{equation}
which represents  an additional fundamental constraint between stress and strain rate explicitly independent of viscosity coefficients.  Trivial rearrangement of \eqref{fundrel} gives
\begin{equation}
\label{actfund}
\frac{\left(\sigma_{11}-\sigma_{22}\right)}{\sigma_{12}}=\frac{\left(e_{11}-e_{22}\right)}{e_{12}}.
\end{equation}

It is well-known that a solution to the system \eqref{sequi}  expressed in terms of the Airy stress function $\chi(x,y)$ is represented by:
\begin{equation}
\label{asf}
\sigma_{11}=-\chi_{,22},\qquad \sigma_{22}=-\chi_{,11},\qquad \sigma_{12}=\chi_{,12}.
\end{equation}

Substitution of  \eqref{asf} in \eqref{actfund} yields 
\begin{equation}
\label{apde}
\chi_{,11}-\chi_{,22}-\Lambda(x_{1},\,x_{2})\chi_{,12}=0,\qquad  \Lambda:=\frac{(e_{11}-e_{22})}{e_{12}},
\end{equation}
which is the partial differential equation satisfied by the Airy stress function  $\chi(x,y)$ explicitly  independent of viscosity coefficients  $\lambda,\, \mu$.

On the assumption that \eqref{apde} can be solved for general $(v_{1},\,v_{2})$ and therefore $\Lambda$, the solution may be used to derive the equilibrium stress components and consequently the viscosity coefficient  $\mu$ from \eqref{fundrel}$_{3}$. Furthermore, it follows  from \eqref{asf} and \eqref{trace} that
\begin{equation}
\label{hydro}
\chi_{,\alpha\alpha}=-2(\lambda+\mu)(e_{11}+e_{22}),
\end{equation}
which may be used to determine $(\lambda+\mu)$ and therefore  $\lambda(x_{1},x_{2})$.
\section{Derivation of viscosity coefficients in Theorem~\ref{nsthm}}\label{A2}
Derivation of expressions \eqref{mu} and \eqref{lamb} for the viscosity coefficients  stipulated in Theorem~\ref{nsthm} is conveniently described
in terms of the complex variable $z$ and its conjugate $\bar{z}$  defined by 
\begin{equation}
\label{complex}
z=(x_{1}+ix_{2}),\qquad \bar{z}=(x_{1}-ix_{2}).
\end{equation}

Details of the following computations are contained in standard texts; e.g., Muskhelishvili \cite{musk63}. 

Differentiation with respect to $z$ and $\bar{z}$ is  defined to be 
\begin{eqnarray}
\label{zdiff}
\frac{\partial}{\partial z}&:=& \frac{1}{2}\left(\frac{\partial}{\partial x}-i\frac{\partial}{\partial y}\right),\\
\label{zbardiff}
\frac{\partial }{\partial\bar{z}}&:=& \frac{1}{2}\left(\frac{\partial}{\partial x}+i\frac{\partial}{\partial y}\right).
\end{eqnarray}


 The particular velocity components \eqref{orb1} and \eqref{orb3}  of present concern, repeated for convenience 
\begin{equation}
\label{vel}
v_{1}=x_{2},\qquad v_{2}=\frac{(x^{2}_{2}-x_{1}^{2})}{2x_{1}},\qquad (x_{1},\,x_{2}) \in \Omega_{1},
\end{equation}
are the real and imaginary parts of   the (non-analytic) velocity field $v(z,\,\bar{z})$   represented by 
\begin{equation}
\label{vexp}
v= v_{1}+iv_{2}=-i\frac{z^{2}}{\left(z+\bar{z}\right)}.
\end{equation}
The corresponding strain rate components become 
\begin{eqnarray}
\label{deone}
e_{11}&=& 0,\\
\label{detwo}
e_{22}&=& i\frac{\left(\bar{z}-z\right)}{\left( z+\bar{z}\right)},\\
\label{de12}
e_{12}&=& \frac{1}{2}\frac{\left(z^{2}+\bar{z}^{2}\right)}{\left(z+\bar{z}\right)^{2}}.
\end{eqnarray}

Insertion into the formula for $\Lambda(z,\bar{z})$ (see \eqref{apde}$_{2}$) gives
\begin{eqnarray}
\label{dlam}
\Lambda= 2i\frac{\left(z^{2}-\bar{z}^{2}\right)}{\left(z^{2}+\bar{z}^{2}\right)}.
\end{eqnarray}
Consequently, \eqref{apde} assumes the form
\begin{equation}
\label{dcomapde}
\left(\chi_{,zz}+\chi_{,\bar{z}\bar{z}}\right)\left(z^{2}+\bar{z}^{2}\right)+\left(\chi_{,zz}-\chi_{,\bar{z}\bar{z}}\right)\left(z^{2}-\bar{z}^{2}\right)=0,
\end{equation}
which upon rearrangement simplifies to
\begin{equation}
\label{findcomapde}
z^{2}\chi_{,zz}+\bar{z}^{2}\chi_{,\bar{z}\bar{z}}=0.
\end{equation}

Among the possible solutions to \eqref{findcomapde}, select that obtained from  setting
\begin{equation}
\label{spsep}
z^{2}\chi_{,zz}=k=-\bar{z}^{2}\chi_{,\bar{z}\bar{z}},
\end{equation} 
for real constant $k$. Integration of the equation on the left gives
\begin{equation}
\label{genspsep}
\chi(z,\bar{z})=-k\log{z}+zf(\bar{z}) +q(\bar{z}),
\end{equation}
where $f(.),\,q(.)$ are arbitrary functions. Similarly, integration of  the equation on the right of \eqref{spsep} yields
\begin{equation}
\label{conspsep}
\chi(z,\bar{z})= k\log{\bar{z}} +\bar{z}g(z)+p(z),
\end{equation}
where $g(.),\,p(.)$ are arbitrary functions.
Choose
\begin{eqnarray}
\nonumber
p(z)&=& -k\log{z},\qquad q(\bar{z})=k\log{\bar{z}},\\
\nonumber
g(z)&=& 2z\log{z},\qquad f(\bar{z})=-2\bar{z}\log{\bar{z}}.
\nonumber
\end{eqnarray}
to derive  respectively from \eqref{genspsep} and \eqref{conspsep} the equivalent expressions
\begin{eqnarray}
\nonumber
\chi(z,\bar{z})&=& -k\log{\frac{z}{\bar{z}}}-2z\bar{z}\log{\bar{z}},\\
\nonumber
\chi(z,\bar{z})&=& -k\log{\frac{z}{\bar{z}}}+2z\bar{z}\log{z},
\end{eqnarray}
and consequently  the final form for the Airy stress function $\chi(z,\bar{z})$  is:
\begin{equation}
\label{finchi}
\chi(z,\bar{z})=\left(z\bar{z}-k\right)\log{\frac{z}{\bar{z}}}.
\end{equation}
Substitution in \eqref{asf} 
determines the equilibrium stress components as
\begin{eqnarray}
\label{dsigone}
\sigma_{11}(z,\bar{z})&=& i\left[(z\bar{z}+k)\frac{(z^{2}-\bar{z}^{2})}{(z\bar{z})^{2}}+2\log{\frac{z}{\bar{z}}}\right],\\
\label{dsigtwo}
\sigma_{22}(z,\bar{z})&=& i\left[(z\bar{z}+k)\frac{(\bar{z}^{2}-z^{2})}{(z\bar{z})^{2}}+2\log{\frac{z}{\bar{z}}}\right],\\
\label{dsig12}
\sigma_{12}(z,\bar{z})&=& \frac{(z\bar{z}+k)(z^{2}+\bar{z}^{2})}{(z\bar{z})^{2}},
\end{eqnarray}
while the viscosity coefficient $\mu(z,\,\bar{z})$ from \eqref{fundrel}$_{3}$ and \eqref{de12} is given by 
\begin{equation}
\label{dmu}
\mu= \frac{(z\bar{z}+k)(z+\bar{z})^{2}}{(z\bar{z})^{2}}.
\end{equation}
In view of \eqref{deone}, \eqref{deone},\eqref{detwo} and \eqref{dsigone}, the viscosity coefficient  $\lambda(z,\,\bar{z})$ is
\begin{eqnarray}
\nonumber
\lambda&=&\frac{\sigma_{11}}{e_{22}}\\
\nonumber
&=& \frac{2(z+\bar{z})}{(\bar{z}-z)}\log{\frac{z}{\bar{z}}}-\frac{(z\bar{z}+k)(z+\bar{z})^{2}}{(z\bar{z})^{2}}.
\end{eqnarray}
In terms of polar coordinates $(r,\theta)$, where
\begin{equation}
\nonumber
z=r\exp{(i\theta)},\qquad \bar{z}=r\exp{(-i\theta)},
\end{equation}
these expressions are written as 
\begin{eqnarray}
\nonumber
\sigma_{11}&=& -2\left[2\theta+\left(\frac{(r^{2}+k)}{r^{2}}\sin{2\theta}\right)\right],\\
\nonumber
\sigma_{22}&=& 2\left[-2\theta+\left(\frac{(r^{2}+k)}{r^{2}}\sin{2\theta}\right)\right],\\
\nonumber
\sigma_{12}&=& 2\left(\frac{(r^{2}+k)}{r^{2}}\right)\cos{2\theta},\\
\nonumber
e_{11}&=& 0,\\
\nonumber
e_{22}&=& \tan{\theta},\\
\nonumber
e_{12}&=& \frac{\cos{2\theta}}{4\cos^{2}{\theta}},\\
\nonumber
\mu&=& 4\frac{(r^{2}+k)}{r^{2}}\cos^{2}{\theta},\\
\nonumber
\lambda&=& -4\theta \cot{\theta}-4\left(\frac{(r^{2}+k)}{r^{2}}\right)\cos^{2}{\theta}.
\end{eqnarray}
It may easily be checked by direct substitution that the stress is in equilibrium under zero body force and that
\begin{equation}
\nonumber
v_{\alpha}n_{\alpha}=0
\end{equation}
at all points on $\partial\Omega_{1}\backslash (0,0)$.

 The proof of \eqref{mu} and \eqref{lamb}  is complete on taking $k=0$ in the above expressions.

\section{Conservative Lagrange trajectories}\label{A3}
Let $(x,y)$ be the rectangular coordinates of a point moving with respect to time and suppose that the differentiable  function $G(x,y)$ is conserved so that $G(x(t),\,y(t))$ is constant.  Examples are the entropy in an adiabatic system or a one parameter family of plane closed curves discussed in Section~\ref{mot}. 
Define velocities $(u,\,v)$ by
\begin{eqnarray}
\label{g1}
u(x,t)&:=&\dot{x},\\
\label{g2}
v(x,t)&:=&\dot{y},
\end{eqnarray}
to obtain from 
\begin{equation}
\nonumber
\frac{d G(x(t),\,y(t))}{dt}=0
\end{equation}
the relation
\begin{equation}
\label{gder}
\frac{\partial G}{\partial x}u+\frac{\partial G}{\partial y}v=0.
\end{equation}
Consequently, \emph{any} $u$   determines $v$ from the expression
\begin{equation}
\label{vdef}
v=-\frac{\partial G}{\partial x}\left(\frac{\partial G}{\partial y}\right)^{-1}u
\end{equation}
such that the system \eqref{g1} and \eqref{g2} is conservative with $G(x,y)$ as first integral.

The particular choice
\begin{equation}
\label{uham}
u(x,y)= \frac{\partial G}{\partial y}
\end{equation}
yields  a Hamiltonian system for which
\begin{equation}
\label{vdef}
v= -\frac{\partial G}{\partial x}.
\end{equation}
It follows that the incompressibity condition
\begin{equation}
\label{incompre}
\frac{\partial u}{\partial x}+\frac{\partial v}{\partial y}=0
\end{equation}
 holds provided the order of the second partial derivatives of $G$ can be reversed.

For fluids with variable mass density, the continuity equation is satisfied  by  modifying the definition of $u$.  Thus, suppose 
\begin{equation}
\label{guw}
u:= \frac{\partial G}{\partial y}w(x,\,y),
\end{equation}
where $w(x,\,y)$ is some sufficiently smooth function, and $G(x,\,y)$ continues to be the first integral of the conservative system \eqref{g1} and \eqref{g2}. Accordingly, $v(x,\,y)$ is modified to
\begin{equation}
\label{gvw}
v:= -\frac{\partial G}{\partial x}w(x,\,y),
\end{equation}
and  the continuity equation 
\begin{equation}
\label{gcont}
\frac{\partial}{\partial x}(\rho u)+\frac{\partial }{\partial y}(\rho v)=0,
\end{equation}
is satisfied provided when $(i)$ the mass density $\rho(x,y)$ is chosen to be  $\rho=D w^{-1}$, for constant $D$, and $(ii)$  the order of  the second partial derivatives of $G$ can be reversed.
\begin{rem}
Apart from differentiability of $G$, no further assumptions have been introduced on either $G(x,\,y)$ or $w(x,\,y)$, whose choice clearly determines the smoothness of the mass density $\rho(x,\,y)$.
\end{rem}
\begin{rem}
Lagrange trajectories  corresponding to  the system \eqref{eu1} and \eqref{eu2} are a special case of \eqref{g1} and \eqref{g2} upon selecting $w(x,y) =x$ and setting
\begin{equation}
\label{spG}
 G(x,\,y)=\frac{(x^{2}+y^{2})}{2x} \quad \left(=H(x,y)\right).
\end{equation}
The mass density becomes $\rho=x^{-1}$ and
\begin{equation}
\label{gspder}
\frac{\partial G}{\partial x}=\frac{(x^{2}-y^{2})}{2x^{2}},\qquad \frac{\partial 
G}{\partial y}=\frac{y}{x}.
\end{equation}
which on insertion into \eqref{guw} and \eqref{gvw} yields the required velocity components. Note that 
\begin{equation}
\nonumber
\frac{\partial ^{2}G}{\partial x\partial y}\neq \frac{\partial ^{2}G}{\partial y\partial x}\qquad \mbox{ for }  (x,\,y) =(0,\,0).
\end{equation}
\end{rem}
\begin{rem}
Substitution of \eqref{guw} and \eqref{gvw} in the Euler equations yields  equations for the partial derivatives of the pressure which can be solved for the pressure.  Consequently, for any conserved quantity $G(x,y)$ and any choice of $w(x,\,y)$ the Euler equations can be solved for appropriate pressure determined semi-inversely. 

When $w=x$ the pressure is given by \eqref{epress} which to within a constant is the  conserved quantity $G(x,y)$ given by \eqref{spG}. 
\end{rem}
\begin{rem}
For the choice $w(x,\,y)=x^{\alpha},\,0<\alpha <1$, the density becomes $\rho =x^{-\alpha}$ and though singular at the origin, is integrable on the region $\Omega$ and computations involving mass density are valid. 
\end{rem}

\end{document}